# Data envelopment analysis models or the virtual gap analysis model: Which should be used for identifying the best benchmark for each unit in a group?

Professor Emeritus Fuh-Hwa Franklin Liu[1]
Department of Industrial Engineering and Management, National Yang Ming Chiao Tung University
No. 1001, Daxue Rd. East Dist., Hsinchu City 30093, Taiwan, Republic of China

Associate Professor, Su-Chuan Shih
Department of Business Administration, Providence University
200, Sec. 7, Taiwan Boulevard, Shalu Dist., Taichung City 43301, Taiwan, Republic of China

---

[1] Corresponding author.
E-mail: fliu@nycu.edu.tw
**Author Contribution Statement**:
Author A: Conceptualization, Methodology, Writing-original draft.
Author B: Conceptualization, Data curation, Validation, Formal Analysis, Visualization, Software.

# Data envelopment analysis models or the virtual gap analysis model: Which should be used for identifying the best benchmark for each unit in a group?

## Abstract

Decision-making units (DMUs) in a group convert the same resources (i.e., input indices) into the same products (i.e., output indices) at different scales. Performance indices have different measurement units, and their market prices per unit are unobtainable. Data envelopment analysis (DEA) programs employ linear programming to estimate the *virtual weight* and best *slack* of every input and output index for each DMU, named DMU-o, to obtain the minimum relative inefficiency against the DMUs. DMU-o reduces each input's slack, the surplus, and expands each output's slack, the shortage, to the benchmark. Each DEA program specifies an *artificial goal weight* for each performance index. The relative inefficiencies in the primal and dual models are the sum of the weighted slacks and the virtual gap of the total virtual weighted inputs to the outputs, respectively. DEA programs have failed the uncountable attempts to conceive the artificial goal weight equal to the estimated virtual weight for each performance index; therefore, they have incomplete solutions that some of the slacks could not be aggregated into the efficiency score. Our new virtual gap analysis program assesses DMU-o comprehensively. The four-phase procedure ensures DMU-o has the achievable best benchmarks for implementation and its compatible best peers to learn. Each DMU is a point in the 2D geometric intuition of the virtual technology set in assessing DMU-o. The best peers and the improved DMU-o are on the best efficiency boundary. Inefficient DMUs are situated underneath the boundary.

## 1. Research background and objectives

### 1.1 Assessing each unit in a group

Decision-making units (DMUs) in a group convert the same resources (i.e., input indices) into the same products (i.e., output indices) at different scales. Performance indices have different measurement units, and their market prices per unit are unobtainable. In the past 45 years, *data envelopment analysis* (DEA) models in tens of thousands of journal articles and hundreds of textbooks have attempted to assess DMUs using the observed datasets of inputs and outputs of DMUs. DEA models are based on neoclassical economics production theory and have incomplete solutions. The production theory and hypotheses are well illustrated in the productivity and efficiency measurement literature, such as in the newest comprehensive textbook by Sickles and Zelenyuk (2019).

We propose a *linear programming* (LP)-based virtual gap analysis (VGA) program to assess each DMU, named DMU-o, against peers in a group. The VGA program satisfies Criteria (i)–(v): (i) Simultaneously estimate the virtual unit price and best benchmark of every input and output index. Use the solutions to compute the maximum relative efficiency that is ensured between 0 and 1. (ii) DMU-o reduces each input's surplus and expands each output's shortage to its estimated benchmark to improve the relative efficiency to 1. (iii) The estimated virtual unit price of each input and output is applied to its current amount and estimated slack; (iv) DMUo increases the current relative efficiency to 1 by changing the scale. Hence, the relative efficiency should include a scaling efficiency to reflect a scaling efficiency; (v) The interval of each input and output index's benchmark is estimated for DMUo to determine the achievable level for implementation. The selected relative efficiency and scale efficiency are applicable.

*1.2 Data envelopment analysis programs: Incomplete assessments*

Charnes, Cooper, and Rhodes (CCR) (1978) introduced a remarkable linear programming (LP)-based DEA program to achieve Criteria (i) and (ii). Banker, Charnes, and Cooper (BCC) (1984) added the convexity condition to the CCR program and yielded the BCC program for measuring scale efficiency, as mentioned in Criterion (iv). The dual variable corresponds to the convexity condition denoted as $w_o$. The solution $w_o^* < 0$, $w_o^* = 0$, and $w_o^* > 0$ indicates that DMU-o exhibits increasing, constant, and decreasing *returns to scale* (RTS). The programs without and with the convexity condition are called the *constant RTS* (CRS) and *variable RTS* (VRS) programs, respectively.

The CCR and BCC programs estimate the radial efficiency that satisfies Criterion (i). Criterion (ii) cannot be fulfilled when some surpluses or shortages are not aggregated into the radial efficiency. Countless non-orientation, non-radial (NONR) DEA programs have attempted to conquer the drawbacks of CCR and BCC programs. Halická and Trnovská (2021) comprehensively reviewed the literature on NONR DEA programs. There are multiple DEA programs for assessing DMU-o. How to choose between them?

DEA programs estimate the *virtual weight* and best *slack* of every input and output index for DMU-o to obtain the relative inefficiency against the efficiency frontier of DMUs. DMU-o reduces each input's slack, the surplus, and expands each output's slack, the shortage, to the benchmark. Each DEA program specifies an *artificial goal weight* for each slack. The relative inefficiency in the primal and dual models is the sum of the weighted slacks and the virtual gap of the total virtual weighted inputs to the outputs. In the dual model, each artificial goal weight is a conditioned lower bound to the virtual weight.

Numerous DEA programs have conceived different artificial goal weights that are unequal to the estimated virtual weights. Further, DEA assumes that the CRS and VRS technology sets have CRS and VRS technological efficiency frontiers, respectively. The same artificial goal weights were used in each DEA program to assess every DMU-o under CRS and VRS conditions. The estimated virtual weights could not equal the artificial goal weights because the DEA programs were ill-formulated based on impractical assumptions. The DEA programs obtained incomplete solutions.

In the CCR and BCC programs, the non-Archimedean infinitesimal $\varepsilon$ represents the artificial goal unit price of surplus and shortage, which is unequal to each input and output's estimated virtual unit price. Halická and Trnovská (2021, Tables 2 and 3) listed the objective functions and lower bounds of virtual unit prices, the goal unit prices, of various NONR DEA programs. They concluded that some VRS NONR DEA programs might provide an efficiency score not in (0, 1). Some VRS NONR programs may not meet some of the five desired properties.

In the CRS NONR DEA programs, DMU-o satisfies Criterion (i) only if the total price of inputs equals 1, and the remaining criteria cannot be satisfied. The slack-based measurement (SBM) program introduced by Tone (2001) could satisfy Criterion (i) only. Our numerical example in Section 5 demonstrates this. DEA programs have incomplete assessments and cannot fulfill Criteria (i)–(iv) in evaluating every DMU-o.

Banker, Charnes, and Cooper (1984) provided a generic picture of the 2D (Figure 2) technology set of DMUs with one input and one output indices. The piecewise segments on the VRS efficient technological frontier portrayed three situations of scale changes: increased RTS, constant RTS, and decreased RTS. The three situations have $w_o^* < 0$, $w_o^* = 0$, and $w_o^* > 0$ for the intercept values associated with the three tangent lines of the points at the three segments. Figure 2 portrays the tangent lines intercepting at the horizontal, origin, and vertical axes. The measurement unit of $w_o^*$ is the same as the input index if $w_o^* < 0$ and is the same as the output index if $w_o^* > 0$. This phenomenon reveals that the 2D illustration is inadequate. The assumption that the observed dataset, the technology set, has a VRS frontier cannot be true. This research defines a 2D geometric intuition of the *virtual technology set* in assessing DMU-o. The *best efficiency boundary* and an anchor point for the estimated $w_o^*$ could be visualized (see Section 5 of the numerical example).

Banker, Charnes, and Cooper (1984, Figure 3) used 2D geometric intuition with one input and one output indices to demonstrate the relationships between the input-oriented CRS and VRS frontiers. There is no envelope for the production possibility set with more than two performance indices. Cooper, Seiford, and Tone (2007, Definition 5.2) derived scale efficiency equals the CCR relative efficiency ratio to the BCC relative efficiency. The CCR and BCC radial efficiencies may not aggregate all slacks for DMUs with multiple inputs and outputs; hence the scale efficiency is inaccurate. Figure 3 fails to illustrate the CRS and VRS frontiers related to scale efficiency.

Similarly, the other CRS and VRS NONR DEA programs could not obtain scale efficiency. Furthermore, DEA theory added the convexity condition to the CRS production possibility set to form the VRS production possibility set. The hyperplanes of the CRS and VRS primal DEA programs, respectively, without and with the convexity condition, construct the convex polyhedral sets based on the theory of linear programming. DEA theory misused the convexity condition to link the CRS and VRS DEA models. At the same time, the dual CRS and VRS programs are convex polyhedral sets that have not to do with the dual variable $w_o$ that is associated with the convexity condition. The decision variables of the production possibility set and polyhedral set are different. We have developed a new method to derive scale efficiency from satisfying Criterion (iv).

The CCR, BCC, and NONR programs are fundamental to countless DEA programs that could not satisfy Criterion (iii). DEA programs are ill-formulated in many characteristics. Emrouznejad and Yang (2017) cited over 10,000 journal articles and over 100 textbooks for DEA-related research. Dyson et al. (2001), Spinks et al. (2009), Brown (2006), Zhu et al. (2018), Zarrin et al. (2022), and many other researchers investigated the pitfalls and protocols of DEA, and Mergoni and De Witte (2022) and Panwar et al. (2022) conducted systematic literature reviews of DEA.

*1.3 Virtual gap analysis programs: Comprehensive assessment of each DMU*

Instead of artificial goal weights, our two new *linear programming* (*LP*)-based VGA programs postulate the "virtual goal unit price" and systematically determine them in two steps. The final virtual goal unit prices equal the estimated virtual unit prices to assess DMU-o comprehensively; hence, Criterion (iii) could be satisfied. The estimated total *virtual price of inputs* and the t*otal virtual price of outputs* of each DMU are denoted as *vInput* and vOutput. The primal model aims to have the minimum virtual gap price of DMU-o, which equals the excesses of vInput to vOutput. We defined the *sum of intensities condition* (SIC) as the sum of DMUs' variable intensities equal to the given alternative scalar, $\kappa_o^a$. The relative inefficiency of the first VGA model is the *pure technical inefficiency* (PTE), which is the virtual gap price ratio to the vInput of DMU-o. Thus, we named the first program the PTE program. The first VGA program without the SIC could not satisfy Criterion (iv) and Criterion (v).

In the dual model of the second VGA program, the SIC equals a given alternative SIC scalar, $\kappa_o^a$. The objective function of the primal model equals the vInput minus the vOuput, plus the additional element, $\kappa_o^a w_o$, which is the *virtual price of the SIC scalar* (*vScalar*). The proportions of the total slack prices of inputs to outputs in the objective function of the dual model are $\gamma_o$ and $(1 - \gamma_o)$. The estimated vScalar is partitioned into two parts, $\gamma_o \kappa_o^a w_o^*$ and $(1 - \gamma_o)\kappa_o^a w_o^*$, to reflect the effects on the input and output indices.

The *affected vInput* (*avInput*) equals the vInput plus $\gamma_o \kappa_o^a w_o^*$, and the *affected vOutput* (*avOutput*) equals the vOutput minus $(1 - \gamma_o)\kappa_o^a w_o^*$. Thus, the virtual gap price of DMU-o becomes the avInput minus the avOutput. Similarly, the avInput and avOuput express each DMU's virtual technology. Similar to the PTE program, the 2D geometric intuition illustrates the solutions of the second VGA program. The anchor point represents vScalar and is located at the coordinate $(\gamma_o \kappa_o^a w_o^*, [1 - \gamma_o]\kappa_o^a w_o^*)$ in the second and fourth quadrants when $w_o^* < 0$ and $w_o^* > 0$, respectively. The relative inefficiency equals the estimated virtual gap price ratio to the avInput. The scale efficiency equals the vScalar ratio to the avInput. The technical inefficiency is the relative inefficiency minus the scale efficiency.

The second VGA program has the SIC, and the relative inefficiency includes *scale and technical inefficiencies* (STE). Thus, we named the second VGA program the *STEa* model with the alternative scalar, $\kappa_o^a$. The 2D geometric intuition depicts the relations of the relative, technical, and scale efficiencies to the vectors between the origin O, anchor, and DMU-o points, respectively.

To satisfy Criterion (v), DMU-o takes a for-phase procedure to determine the scalar, $\kappa_o^a$. The four-phase procedure comprises the PTE, STE1, STE2, and STE^ programs. In Phase 1, the PTE program provides the first SIC scalar, $\kappa_o^1$. In Phase 2, the STE1 program has the SIC with scalar $\kappa_o^1$. The STE1 program quantifies the scale efficiency of the PTE program. Phase 3 identifies the second SIC scalar, $\kappa_o^2$, according to the standard sensitivity analysis on the right-hand side scalar $\kappa_o^1$ of the STE1 program. In Phase 4, each best peer and DMU-o have compatible production technologies, as confirmed by DMU-o. The 2D geometric intuition helps DMU-o to select the target SIC scalar, $\hat{\kappa}_o$, between the $\kappa_o^1$ and $\kappa_o^2$. The STE^ program should provide an achievable benchmark for each input and output for implementation. Modifying the VGA programs and procedures could satisfy further considerations in research and applications.

Using the estimated virtual unit prices could compute each DMU's avInput and avOutput. Each DMU is a point at the coordinate (avInput, avOutput) in the 2D geometric intuition of the virtual technology set. In the 2D geometric intuition, the diagonal at the origin O is the best efficiency boundary. Points of benchmark peers are located on the boundary, and the points of inefficient DMUs are situated underneath the boundary. DMU-o projects on the best efficiency boundary, and its relative boundary efficiency equals 1. The 2D geometric intuition depicts that relative efficiency equals technical efficiency minus scale efficiency.

### *1.4 Our early VGA developments*

Hwang (2014) worked with the first author to establish the CRS and VRS VGA programs, in which the virtual goal price of *DMU-o* was incorrect. The relationship between the CRS and VRS VGA programs was incorrect. Previous studies have used the CRS VGA programs for production systems with two phases (Liu and Liu, 2017a, 2017b). The first author has advised on several unpublished theses using the CRS and VRS VGA programs for analyzing real datasets. Li (2022) worked with the second author and employed the PTE program to evaluate the hotel industry's internet word-of-mouth performance through online travel reviews. The CRS and VRS VGA programs were modified into the PTE and STEa programs without assumptions and axioms of the production theory. The STE^ program satisfies all five criteria mentioned above.

*1.5 Paper organization*

The remainder of this paper is organized as follows. Section 2 introduces the PTE and STEa programs and their two-step process to solve the virtual goal unit price of *DMU-o*. The post-analysis of the assessment items is presented in Section 3. In Section 4, the four-phase VGA procedure is described. Section 5 offers a typical numerical example of DMUs with multiple inputs and outputs. A 2D geometric intuition could visualize the assessment results. Finally, Section 6 concludes and discusses this paper.

## 2. Virtual gap analysis (VGA) procedure

Sections 2.2 and 2.3 present the PTE and STEa programs, in which the primal and dual models are the *total gap price* (TGP) model and the *total slack price* (TSP) model. The dualities of the VGA programs are discussed in Section 2.4.

*2.1 Notations*

$J$, $I$, and $R$ denote the sets of DMUs, inputs, and outputs, respectively, where $J = \{j|j=1,...,n\}$, $I = \{i|i=1,...,m\}$, and $R = \{r|r=1,...,s\}$. The symbols $x_{ij}$ and $y_{rj}$ represent the non-negative continuously disproportionate *input-i* and *output-r* of *DMU-j*, respectively. The size of the observed dataset ($X$, $Y$) is $m \times n$, in which $x_j$ and $y_j$ are the column vectors. The value of $n$ should be larger than $m + s$.

The decision variables are linked to *input-i* and *output-r* of *DMU-o* (i.e., $v_{io}$ and $u_{ro}$), also called the *virtual unit prices* or weights, in the VGA program symbols $v_o$ and $u_o$ denote the row vectors. The symbols $v_{io}x_{io}$ and $u_{ro}y_{ro}$ represent the virtual prices measured in the *virtual currency* (\$). Both $q_{io}$ and $p_{ro}$ denote the *slack* of *input-i* and *output-r*, respectively, where $Q_{io}$ and $P_{ro}$ are the *slack ratios*, $q_{io}/x_{io}$ and $p_{ro}/y_{ro}$, respectively. Here, $q_o$, $Q_o$, $p_o$, and $P_o$ denote the vector columns; $\pi_j$ denotes the *intensity* of the *DMU-j* for evaluating *DMU-o*, and $\pi$ denotes the vector column.

The symbols $g_{io}^x$ and $g_{ro}^y$ denote the virtual goal unit price of input-i and output-r in assessing DMU-o. The conditions in the TGP model, $v_{io} \geq g_{io}^x$ and $u_{ro} \geq g_{ro}^y$, are equivalent to $x_{io}v_{io} \geq \tau_o$ and $y_{ro}u_{ro} \geq \tau_o$, respectively, where $g_{io}^x x_{io} = \tau_o$ and $g_{ro}^y y_{ro} = \tau_o$, respectively. The symbol $\tau_o$ denotes the systematically determined virtual goal price of DMU-o, which is the normalized lower bound of all the virtual prices of the input and output indices. The virtual slack prices of input-i and output-r are $g^x_{io} q_{io}$ and $g^y_{ro} p_{ro}$ in the TSP model, respectively, and are equivalent to $\tau_o Q_{io}$ and $\tau_o P_{ro}$, respectively. The symbol e denotes the row or column unit vector. The superscripts # and * indicate the Step I and Step II optimal solutions for a decision variable in the VGA programs.

## 2.2 PTE program and the solving process

The PTE program and the particular two-step solving process are presented here.

TSP model of the PTE program: (1)

$$\delta_o^{PTE*} = Max \sum_{i \in I} Q_{io} \tau_o + \sum_{r \in R} P_{ro} \tau_o \quad (1.0)$$

$$s.t. \sum_{j \in J} x_{ij} \pi_j + Q_{io} x_{io} = x_{io}, \forall i \in I; \quad (1.1)$$

$$-\sum_{j \in J} y_{rj} \pi_j + P_{ro} y_{ro} = -y_{ro}, \forall r \in R; \quad (1.2)$$

$$\pi, Q_o, P_o \geq 0. \quad (1.3)$$

TGP model of the PTE program: (2)

$$\Delta_o^{PTE*} = Min \sum_{i \in I} v_{io} x_{io} - \sum_{r \in R} u_{ro} y_{ro} \quad (2.0)$$

$$s.t. \sum_{i \in I} v_{io} x_{ij} - \sum_{r \in R} u_{ro} y_{rj} \geq 0, \forall j \in J; \quad (2.1)$$

$$x_{io} v_{io} \geq \tau_o, \forall i \in I; \quad (2.2)$$

$$y_{ro} u_{ro} \geq \tau_o, \forall r \in R; \quad (2.3)$$

$$v_o, u_o \text{ free}. \quad (2.4)$$

Here, Eq. (2.1) limits the total virtual gap price of each *DMU-j* to be nonnegative, and Eqs. (2.2) and (2.3) restrict each virtual price of *DMU-o* to the lower bound $\tau_o$. In Step I, substitute $\tau_o$ by $\tau_o^\# = 1$, will obtain the optimal solutions of the TSP model (1), $(Q_o^*, P_o^*, \pi^*) = (Q_o^\#, P_o^\#, \pi^\#)$. In Step II, substitute $\tau_o$ by $\tau_o^* = \$ t$. Eqs. (3) and (4) are the TSP and TGP models (1) and (2) in Steps I and II. The TGP equals the *total virtual price of inputs* (vInput), $\alpha_o^\#$ and $\alpha_o^*$, excess to the *total virtual price of outputs* (vOutput), $\beta_o^\#$ and $\beta_o^*$. At the same time, TSP equals $\sum_{i \in I} Q_{io}^* \times \tau_o^*$ plus $\sum_{r \in R} P_{ro}^* \times \tau_o^*$.

$$\Delta_o^{PTE\#} = v_o^\# x_o - u_o^\# y_o = \alpha_o^\# - \beta_o^\# > 0 \text{ and } \delta_o^{PTE\#} = \sum_{i \in I} Q_{io}^\# \times \$1 + \sum_{r \in R} P_{ro}^\# \times \$1 > 0. \quad (3)$$

$$\Delta_o^{PTE*} = v_o^* x_o - u_o^* y_o = \alpha_o^* - \beta_o^* > 0 \text{ and } \delta_o^{PTE*} = \sum_{i \in I} Q_{io}^* \times \tau_o^* + \sum_{r \in R} P_{ro}^* \times \tau_o^* > 0. \quad (4)$$

Determine the dimensionless value of $t$ to ensure $\Delta_o^{PTE*}$ between $0 and $1 and $\alpha_o^* = \$1$. Positive $\Delta_o^{PTE\#}$ and $\Delta_o^{PTE*}$ only if inefficient DMUo. According to $\tau_o^\# : \tau_o^* = \$1 : \$t$ and Eq. (3), yields $\alpha_o^\# : \alpha_o^* = 1 : t$. Use Eq. (5) to obtain the dimensionless value of $t$:

$$t = \$1/\alpha_o^\# = \$1/v_o^\# x_o \text{ and } \tau_o^* = \$t. \quad (5)$$

The solutions of Step II are $v_o^* x_o = 1$ and $u_o^* y_o < 1$. Then, use Eq. (6) to compute the minimized relative inefficiency and maximized relative efficiency within zero and one, $0 \leq F_o^{PTE*} < 1$ and $0 < E_o^{PTE*} \leq 1$. That satisfies Criterion (i).

$$F_o^{PTE*} = \Delta_o^{PTE*}/\alpha_o^* = (\alpha_o^* - \beta_o^*)/\alpha_o^* \text{ and } E_o^{PTE*} = 1 - F_o^{PTE*} = \beta_o^*/\alpha_o^*. \quad (6)$$

Hence, normalizing Step I solutions would have Step II solutions, as shown in Eq. (7).

$$\Delta_o^{PTE*} = t\Delta_o^{PTE\#}, \delta_o^{PTE*} = t\delta_o^{PTE\#}, (v_o^*, u_o^*) = t(v_o^\#, u_o^\#). \quad (7)$$

Eq. (8) is also used to calculate the target of each input-$i$ and output-$r$, $\hat{x}_{io}$ and $\hat{y}_{ro}$:

$$\hat{x}_{io}^{PTE} = \sum_{j \in B_o^{PTE}} x_{ij} \pi_j^* = x_{io}(1 - Q_{io}^*), \forall i \in I; \text{ and } \hat{y}_{ro}^{PTE} = \sum_{j \in B_o^{PTE}} y_{rj} \pi_j^* = y_{ro}(1 + P_{ro}^*), \forall r \in R. \quad (8)$$

Criterion (ii) is the *relative boundary efficiency* should be 1, $bE_o^{PTE*} = \sum_{r \in R} \hat{y}_{ro}^{PTE} u_{ro}^* / \sum_{i \in I} \hat{x}_{io}^{PTE} v_{io}^* = 1$. The symbol $B_o^{PTE}$ denotes the best peers (hereafter BPs) identified in the PTE program. *DMU-o* learns virtual technologies of best peers with the estimated intensity, $\pi_j^*$.

## 2.3 STEa program and the solving process

The *sum of intensities condition* (SIC) is the sum of DMUs' variable intensities equals an alternative scalar $\kappa_o^a$, $\sum_{j \in J} \pi_j = \kappa_o^a$. Add the SIC to the PTE program to yield the STEa program. The SIC corresponds to the free-in-sign dual decision variable $w_o$, the virtual unit price of the SIC *scalar* $\kappa_o^a$. The two-step solving process to determine the virtual goal price, $\tau_o^*$, is presented here.

TSP model of the STEa program: (9)

$$\delta_o^{STEa*} = Max \sum_{i \in I} Q_{io} \tau_o + \sum_{r \in R} P_{ro} \tau_o \quad (9.0)$$

$$s.t. \sum_{j \in J} x_{ij} \pi_j + Q_{io} x_{io} = x_{io}, \forall i \in I; \quad (9.1)$$

$$-\sum_{j\in J} y_{rj}\pi_j + P_{ro}y_{ro} = -y_{ro}, \forall r\in R; \quad (9.2)$$

$$\sum_{j\in J}\pi_j = \kappa_o^a; \quad (9.3)$$

$$\pi, Q_o, P_o \geq 0. \quad (9.4)$$

TGP model of the STEa program: (10)

$$\Delta_o^{STEa*} = Min \sum_{i\in I} v_{io}x_{io} - \sum_{r\in R} u_{ro}y_{ro} + \kappa_o^a w_o \quad (10.0)$$

$$s.t. \sum_{i\in I} v_{io}x_{ij} - \sum_{r\in R} u_{ro}y_{rj} + 1w_o \geq 0, \forall j\in J; \quad (10.1)$$

$$x_{io}v_{io} \geq \tau_o, \forall i\in I; \quad (10.2)$$

$$y_{ro}u_{ro} \geq \tau_o, \forall r\in R; \quad (10.3)$$

$$v_o, u_o, w_o \text{ free}. \quad (10.4)$$

In Steps I and II, substitute $\tau_o$ by $\tau_o^\# = \$1$ and $\tau_o^* = \$t$ will solve TSP model (9) comprehensively $(Q_o^*, P_o^*, \pi^*)=(Q_o^\#, P_o^\#, \pi^\#)$. The symbols $\omega_o^{a\#}$ and $\omega_o^{a*}$ are equal to $\kappa_o^a w_o^\#$ and $\kappa_o^a w_o^*$ representing the *virtual price for given the scalar* (vScalar) in Step I and II. To reflect the effects of the SIC on vOutput and vInput, we partition $\omega_o^{a\#}$ into $\gamma_o\omega_o^{a\#}$ and $(1-\gamma_o)\omega_o^{a\#}$. As depicted in Eq. (11), $\gamma_o$ and $(1-\gamma_o)$ are proportions of the total slack price of inputs and outputs in $\delta_o^{STEa\#}$, respectively, with Step II having the same calculation. If $\omega_o^{a\#}$ is equal to 0, then let $\gamma_o$ be equal to 0.5. Eq. (10.0) has solutions in Steps I and II, as shown in Eq. (12) and (13).

$$\gamma_o : (1-\gamma_o) = \tau_o^\# \sum_{i\in I} Q_{io}^\# : \tau_o^\# \sum_{r\in R} P_{ro}^\# = \tau_o^* \sum_{i\in I} Q_{io}^* : \tau_o^* \sum_{r\in R} P_{ro}^*. \quad (11)$$

$$\Delta_o^{STEa\#} = v_o^\# x_o - u_o^\# y_o + \omega_o^{a\#} = [v_o^\# x_o + (1-\gamma_o)\omega_o^{a\#}] - (u_o^\# y_o - \gamma_o\omega_o^{a\#}) = \alpha_{\omega o}^{a\#} - \beta_{\omega o}^{a\#} = \Delta_o^{STEa\#}$$

$$and\ \delta_o^{STEa\#} = \sum_{i\in I} Q_{io}^\# \times \$1 + \sum_{r\in R} P_{ro}^\# \times \$1. \quad (12)$$

$$\Delta_o^{STEa*} = v_o^* x_o - u_o^* y_o + \omega_o^{a*} = [v_o^* x_o + (1-\gamma_o)\omega_o^{a*}] - (u_o^* y_o - \gamma_o\omega_o^{a*}) = \alpha_{\omega o}^{a*} - \beta_{\omega o}^{a*}$$

$$and\ \delta_o^{STEa*} = \sum_{i\in I} Q_{io}^* \times \$t + \sum_{r\in R} P_{ro}^* \times \$t. \quad (13)$$

$\Delta_o^{STEa\#}$ and $\Delta_o^{STEa*}$ are equal to $\alpha_{\omega o}^{a\#}-\beta_{\omega o}^{a\#}$ and $\alpha_{\omega o}^{a*}-\beta_{\omega o}^{a*}$, respectively. They are the *affected vInput* (avInput) minus *affected vOuput* (avOuput) of Step I and Step II.

Determine the dimensionless value of $t$ to have $\Delta_o^{STEa*}$ between 0 $ and $ 1 and $\alpha_{\omega o}^{a*}$ = $ 1. Then, using Eq. (14) to obtain the minimized relative inefficiency and maximized relative efficiency, $0 \leq F_o^{STEa*} < 1$ and $0 < E_o^{STEa*} \leq 1$. That satisfies Criterion (i).

$$F_o^{STEa*} = \Delta_o^{STEa*} / \alpha_{\omega o}^{a*} = (\alpha_{\omega o}^{a*} - \beta_{\omega o}^{a*})/\alpha_{\omega o}^{a*} \text{ and } E_o^{STEa*} = 1 - F_o^{STEa*} = \beta_{\omega o}^{a*}/\alpha_{\omega o}^{a*}. \tag{14}$$

Because $\tau_o^{\#}: \tau_o^{*} = \$1 : \$t$, therefore, if $\alpha_{\omega o}^{a\#}$ equals $t$, then $\alpha_{\omega o}^{a*}$ equals $1. Use Eq. (15) to obtain the dimensionless value of $t$:

$$t = \$ 1 / \alpha_{\omega o}^{a\#} = \$1 / [v_o^{\#} x_o + (1 - \gamma_o)\omega_o^{a\#}] \text{ and } \tau_o^{*} = \$ t. \tag{15}$$

Normalizing the solutions of Step I by the dimensionless value $t$ would obtain the solutions of Step II, as shown in Eq. (16). The solutions of Step II meet Criterion (iii).

$$(\Delta_o^{STEa*}, \delta_o^{STEa*}, v_o^{*}, u_o^{*}, w_o^{a*}) = t \, (\Delta_o^{STEa\#}, \delta_o^{STEa\#}, v_o^{\#}, u_o^{\#}, w_o^{a*}). \tag{16}$$

Eq. (8) calculates the target scale of *input-i* and *output-r*: $\hat{x}_{io}^{STEa}, \forall i$ and $\hat{y}_{ro}^{STEa}, \forall r$. The STE*a* aggregates all slacks to the relative inefficiency. The Criterion (ii), the relative boundary efficiency should be 1, $bE_o^{STEa*} = \sum_{r \in R} \hat{y}_{ro}^{STEa} u_{ro}^{*} / \sum_{i \in I} \hat{x}_{io}^{STEa} v_{io}^{*} = 1$. The symbol $B_o^{STEa}$ denotes the set of best peers that $\pi_j^{*} > 0$. Section 4 illustrates the method to determine the target SIC scalar.

## 2.4 Duality properties

Dantzig and Thapa (1997) proposed various sensitivity analyses to identify the effect of coefficient changes on the optimal solution in a linear program. These sensitivity analyses are crucial when the observed dataset of DMUs is not precisely known. The primal model is the minimized TGP in Eq. (2.0) and (10.0), $\Delta_o^{PTE*}$ and $\Delta_o^{STEa*}$, respectively—the maximized TSP in Eqs. (1.0) and (9.0) are $\delta_o^{PTE*}$ and $\delta_o^{STEa*}$, respectively. Eq. (17) shows that the solutions of the primal and dual objectives are equal in the PTE and STEa programs. The TGVP and TSP are converted into relative inefficiency by Eqs. (6) and (14) for the PTE and STEa programs, respectively.

$$\delta_o^{PTE*} = \Delta_o^{PTE*} \text{ and } \delta_o^{STEa*} = \Delta_o^{STEa*}. \tag{17}$$

The pair of Eqs. (18) and (19) are equivalent to the pair of Eqs. (1.1) and (1.2) or of Eqs. (9.1) and (9.2).

$$\sum_{j \in J} x_{ij} \pi_j = (1 - Q_{io}) x_{io}, \forall i \in I. \tag{18}$$

$$\sum_{j \in J} y_{rj} \pi_j = (1 + P_{ro}) y_{ro}, \forall r \in R. \tag{19}$$

The left-hand sides of Eqs. (18) and (19) form the best efficiency boundary of *DMU-o*. The right-hand sides are the targets of the input and output indices. The equations ensure that all possible slacks aggregate into the objective function, $\sum_{i \in I} Q_{io}^* \tau_o^* + \sum_{r \in R} P_{ro}^* \tau_o^*$, where $Q_{io}^*$ and $P_{ro}^*$ are equally weighted by $\tau_o^*$, the virtual goal price of *DMU-o*.

We listed the strong complementary slackness conditions of the VGA programs in the following equations.

$$[\sum_{j \in B_o^{VGA}} x_{ij} \pi_j^* - x_{io}(1-Q_{io}^*)] v_{io}^* = 0, \forall i \in I \text{ and } [\sum_{j \in B_o^{VGA}} y_{rj} \pi_j^* - y_{ro}(1+P_{ro}^*)] u_{ro}^* = 0, \forall r \in R. \tag{20}$$

In Eqs. (2.4) and (10.4), $v_o$ and $u_o$ are free-in-signs. The non-negative dataset (*X*, *Y*) and $\tau_o^* > 0$ yields the positive $v_o^* > 0$ and $u_o^* > 0$. Therefore, $[\sum_{j \in B_o^{VGA}} x_{ij} \pi_j^* - x_{io}(1 - Q_{io}^*)] = 0, \forall i \in I$ and $[\sum_{j \in B_o^{VGA}} y_{rj} \pi_j^* - y_{ro}(1+P_{ro}^*)] = 0, \forall r \in R$ in set $B_o^{VGA}$, which represents $B_o^{PTE}$ or $B_o^{STEa}$.

$$(v_o^* x_j - u_o^* y_j) \pi_j^* = 0 \text{ and } (v_o^* x_j - u_o^* y_j + 1 w_o^*) \pi_j^* = 0, \forall j \in J. \tag{21}$$

When $(v_o^* x_j - u_o^* y_j) = 0$, *DMU-j* belongs to the best peers of *DMU-o*, $B_o^{PTE} = \{j | \pi_j^* > 0, \forall j\}$. The symbol $w_o^*$ is the unit price of the SIC scalar, and $1 \times w_o^*$ is the *vScalar* of *DMU-j*. When $(v_o^* x_j - u_o^* y_j + 1 w_o^*) = 0$, *DMU-j* belongs to $B_o^{STEa} = \{j | \pi_j^* > 0, \forall j\}$.

$$(v_{io}^* x_{io} - \tau_o^*) Q_{io}^* = 0, \forall i \in I \text{ and } (u_{ro}^* y_{ro} - \tau_o^*) P_{ro}^* = 0, \forall r \in R. \tag{22}$$

When $Q_{io}^* = 0$ and $P_{ro}^* = 0$, then $(v_{io}^* x_{io} - \tau_o^*) > 0$ and $(u_{ro}^* y_{ro} - \tau_o^*) > 0$. In contrast, when $v_{io}^* x_{io} = \tau_o^*$ and $u_{ro}^* y_{ro} = \tau_o^*$, $Q_{io}^* > 0$ and $P_{ro}^* > 0$, respectively.

$$(\sum_{j \in B_o^{STEa}} \pi_j^* - \kappa_o^a) w_o^* = 0. \tag{23}$$

Since $(\sum_{j \in B_o^{STEa}} \pi_j^* - \kappa_o^a) = 0$, *DMU-o* exhibits decreasing, increasing, and constant RTS and the best returns to practice (bRTP) when $w_o^* > 0$, $w_o^* < 0$, and $w_o^* = 0$, respectively.

## 3. Post-analysis of the assessment items

Given the SIC scalar, $\kappa_o^a$, in the STEa program, the solutions can quantify the assessment items. The four-phase procedure illustrated in the subsequent section helps *DMU-o* to select the target scalar.

## 3.1 Virtual technology set in each assessment

Eq. (24) shows the solution for each *DMU-j* in Eq. (10.1). Similar to Eq. (12), which divides $1w_o^*$ into two parts, we rewrite Eq. (24) as Eq. (25).

$$\Delta_j^{STEa*} = (v_o^* x_j - u_o^* y_j) + 1w_o^* \geq 0, \ \forall j \in J\text{-}\{o\}. \tag{24}$$

$$\Delta_j^{STEa*} = [v_o^* x_j + 1 \times (1 - \gamma_o)w_o^*] - (u_o^* y_j - 1 \times \gamma_o w_o^*) = \alpha_{\omega j}^{a*} - \beta_{\omega j}^{a*} \geq 0, \ \forall j \in J\text{-}\{o\}. \tag{25}$$

We denoted the virtual technology sets in the PTE and STEa programs for assessing *DMU-o* as $\Gamma_o^{PTE}$ and $\Gamma_o^{STEa}$, respectively, as shown in Eq. (26), where $(\alpha_j^*, \beta_j^*)$ is equivalent to $(v_o^* x_j, u_o^* y_j)$, and $(\alpha_{\omega o}^{a*}, \beta_{\omega o}^{a*})$ is shown in Eq. (12).

$$\Gamma_o^{PTE} = \{(\alpha_j^*, \beta_j^*), \ \forall j \in J\} \text{ and } \Gamma_o^{STEa} = \{(\alpha_{\omega j}^{a*}, \beta_{\omega j}^{a*}), \ \forall j \in J\}. \tag{26}$$

## 3.2 Best return to practice

*DMU-o* reduced the inputs and expanded the outputs at different rates to the target at the best efficiency boundary. Given the virtual unit prices of the input and output indices, multiple inputs and outputs are aggregated into virtual input and virtual output pairs. The $(\alpha_o^*, \beta_o^*)$ in Eq. (4) represents the virtual technology of *DMU-o* in the PTE program, and $(\alpha_{\omega o}^{a*}, \beta_{\omega o}^{a*})$ in Eq. (13) represent the virtual technology in the STE*a* programs. Similarly, the virtual target technology on the PTE efficiency boundary is represented by $(\hat{\alpha}_o^*, \hat{\beta}_o^*)$ [see Eq. (27)], and the virtual target technology on the STEa efficiency boundary is represented by $(\hat{\alpha}_{\omega o}^{a*}, \hat{\beta}_{\omega o}^{a*})$ [see Eq. (28)].

$$\hat{\alpha}_o^*(\$) = \sum_{i \in I} \hat{x}_{io} v_{io}^* \text{ and } \hat{\beta}_o^*(\$) = \sum_{r \in R} \hat{y}_{ro} u_{ro}^* \tag{27}$$

$$\hat{\alpha}_{\omega o}^{a*}(\$) = [v_o^* \hat{x}_o + (1 - \gamma_o)\omega_o^{a*}] \text{ and } \hat{\beta}_{\omega o}^{a*}(\$) = (u_o^* \hat{y}_o - \gamma_o \omega_o^{a*}) \tag{28}$$

The bRTP in the PTE program is $\Xi_o^{PTE*}$ and in the STEa program is $\Xi_o^{STEa*}$, computed by the following equations.

$$\Xi_o^{PTE*} = (\hat{\beta}_o^*/\beta_o^*) / (\hat{\alpha}_o^*/\alpha_o^*) > 0 \tag{29}$$

$$\Xi_o^{STEa*} = (\hat{\beta}_{\omega o}^{a*}/\beta_{\omega o}^{a*}) / (\hat{\alpha}_{\omega o}^{a*}/\alpha_{\omega o}^{a*}) > 0 \tag{30}$$

The symbol $\Xi_o^{STEa*}$ is the monotone decreasing and increasing of the SIC scalar $\kappa_o^a$ when $w_o^* < 0$ and $w_o^* > 0$, respectively. An efficient *DMU-o* has a bRTP equal to 1.

### 3.3 The virtual technology efficiency boundary in each assessment

Each *DMU-j* in the *virtual technology set* represents a point in the 2D geometric intuition. Figures 1 and 2 illustrate the geometric intuition of the numerical example in assessing the *DMU-o*. The virtual technology set's *best efficiency boundary* is the diagonal at the origin O of the geometric intuition. The best peers are located on the boundary, where *DMU-j* is in the sets of $B_o^{PTE}$ and $B_o^{STEa}$, with $\alpha_j^* = \beta_j^*$ and $\alpha_{\omega j}^{a*} = \beta_{\omega j}^{a*}$, respectively. Simultaneously, the other *DMU-j* lies underneath the boundary, $\alpha_j^* > \beta_j^*$ and $\alpha_{\omega j}^{a*} > \beta_{\omega j}^{a*}$. *DMU-o* projects on the boundary, $\hat{\beta}_o^* = \hat{\alpha}_o^*$ and $\hat{\beta}_{\omega o}^{a*} = \hat{\alpha}_{\omega o}^{a*}$.

The total virtual slack price is the rectilinear distance of *DMU-o* to its projection point on the boundary. The slope of the line segment between the origin O and *DMU-o* is the relative efficiency. The STEa program has an *anchor point* (**APa**) located at the point $([1 - \gamma_o]\kappa_o^a w_o^*, -\gamma_o \kappa_o^a w_o^*)$ in the second and fourth quadrants when $w_o^* > 0$ and $w_o^* < 0$, respectively. The rectilinear distance between APa and the origin O is $\kappa_o^a w_o^*$. The relative, technical, and scale efficiencies can be visualized in the 2D geometric intuition.

### 3.4 Relative, technical, and scale efficiencies

Eqs. (12) and (13) have equivalent expressions in Eq. (31), where $\Delta_{\omega o}^{STEa\#}$ and $\Delta_{\omega o}^{STEa*}$ are the affected virtual technical *gap* price in the Steps I and II evaluations, respectively. The VGP is the sum of the affected technical gap price and *vScalar*, $\Delta_{\omega o}^{STEa*} + \omega_o^{a*}$.

$$\Delta_o^{STEa\#} = (v_o^\# x_o - u_o^\# y_o) + \omega_o^{a\#} = \Delta_{\omega o}^{STEa\#} + \omega_o^{a\#} \text{ and } \Delta_o^{STEa*} = (v_o^* x_o - u_o^* y_o) + \omega_o^{a*} = \Delta_{\omega o}^{STEa*} + \omega_o^{a*}. \quad (31)$$

We divided Eq. (31) by $\alpha_{\omega o}^{a*}$, equivalent to \$1, to become Eq. (32); the relative inefficiency is the affected technical inefficiency plus the scale efficiency. Notably, in Eq. (33), $\ddot{T}_o^{STEa*}$ and $\dot{T}_o^{STEa*}$ are the $(v_o^* x_o - u_o^* y_o)$ divided by $\alpha_{\omega o}^{a*}$ and $v_o^* x_o$, respectively. $\ddot{T}_o^{STEa*}$ may not be between 0 and 1 and $\dot{T}_o^{STEa*}$ is between 0 and 1. Eq. (34) depicts Criterion (iv), the relative efficiency as the affected technical efficiency minus the scale efficiency, $T_o^{STEa*} - S_o^{STEa*}$, and $0 < F_o^{STEa*} = 1 - E_o^{STEa*} \leq 1$.

$$\Delta_o^{STEa*} / \alpha_{\omega o}^{a*} = \Delta_{\omega o}^{STEa*} / \alpha_{\omega o}^{a*} + \omega_o^{a*} / \alpha_{\omega o}^{a*} \Leftrightarrow F_o^{STEa*} = \ddot{T}_o^{STEa*} + S_o^{STEa*}. \quad (32)$$

$$\ddot{T}_o^{STEa*} = (v_o^* x_o - u_o^* y_o) / \alpha_{\omega o}^{a*} \text{ and } \dot{T}_o^{STEa*} = (v_o^* x_o - u_o^* y_o) / v_o^* x_o. \quad (33)$$

$$E_o^{STEa*} = 1 - F_o^{STEa*} = 1 - \ddot{T}_o^{STEa*} - S_o^{STEa*} = T_o^{STEa*} - S_o^{STEa*}. \quad (34)$$

In Section 5, Figure 2 presents the 2D geometric intuition for a numerical example's solutions, where *DMU-o* is DMU-K, and STEa is the STE1 model. The coordinates of points O, K1, and AP1 are (0,0), $(\alpha_{\omega o}^{1*}, \beta_{\omega o}^{1*})$, and $([1 - \gamma_o]\omega_o^{1*}, -\gamma_o \omega_o^{1*})$, respectively. The three rectilinear distances between the points (O,

K1), (AP1, K1), and (AP1, O) are equal to $\Delta_o^{STE1*}$, $\Delta_{\omega o}^{STE1*}$, and $\omega_o^{1*}$, respectively. Figure 2 shows that the sum of the three graphically geometrical vectors is equivalent to $\Delta_o^{STEa*} - \omega_o^{a*} = (v_o^* x_o - u_o^* y_o)$, [see Eq. (31)]. The additions of the three vectors in the horizontal and vertical axes are shown in Eq. (35).

$$\overrightarrow{O,K1} - \overrightarrow{AP1,O} = \overrightarrow{AP1,K1} \Leftrightarrow \alpha_{\omega o}^{1*} - [1-\gamma_o]\omega_o^{1*} = v_o^* x_o \text{ and } \beta_{\omega o}^{1*} + \gamma_o \omega_o^{1*} = u_o^* y_o \qquad (35)$$

The relative, affected technical, and scale efficiency scores are monotone increasing and decreasing from the SIC scalar $\kappa_o^a$ when the solution $w_o^*$ is less than or exceeds 0.

### 3.5 Interlinkage relationships between input and output indices

We rewrote Eq. (12) as Eq. (36) to express the affected virtual prices of *input-i*, $(v_{io}^* x_{io} + \gamma_{io}^Q \omega_o^{a*})$, and *output-r*, $(u_{ro}^* y_{ro} - \gamma_{ro}^P \omega_o^{a*})$. Eq. (37) depicts the proportions of *input-i* and *output-r*, $\gamma_{io}^Q$ and $\gamma_{ro}^P$, respectively, to *vScalar*.

$$\Delta_o^{STEa*}(\$) = [v_o^* x_o + (1-\gamma_o)\omega_o^{a*}] - (u_o^* y_o - \gamma_o \omega_o^{a*}) = \sum_{i \in I}(v_{io}^* x_{io} + \gamma_{io}^Q \omega_o^{a*}) - \sum_{r \in R}(u_{ro}^* y_{ro} - \gamma_{ro}^P \omega_o^{a*}) \qquad (36)$$

$$\gamma_{io}^Q = (1-\gamma_o)Q_{io}^* / \sum_{i \in I} Q_{io}^*, \forall i \in I; \text{ and } \gamma_{ro}^P = \gamma_o P_{ro}^* / \sum_{r \in R} P_{ro}^*, \forall r \in R. \qquad (37)$$

The proportions of $(v_{io}^* x_{io} + \gamma_{io}^Q \omega_o^{a*})$ and $(u_{ro}^* y_{ro} - \gamma_{ro}^P \omega_o^{a*})$ denote the interlinkage relationships between the input and output indices.

## 4. Four-phase VGA procedure

The assessment of *DMU-o* requires four phases that employ the PTE, STE1, STE2, and STE^ programs. Phase–1 determines $\kappa_o^1$, and Phase–2, 3, and 4 use the SIC scalers, $\kappa_o^1$, $\kappa_o^2$, and $\kappa_o^{\hat{}}$, respectively.

### 4.1 Phase 1: Identifying the first scalar, $\kappa_o^1$

Without limiting the DMUs' intensities, the PTE program has the maximized total of the slack ratios in the TSP model [see the model (1)]. *DMU-o* may not achieve the targets computed via Eq. (8). The TGP model [see the model (2)] provides the PTE that could not reflect the characteristics of *DMU-o*'s technology. *DMU-o* should determine its technical and scale efficiencies regarding the best peers' intensities. The PTE program provides an initial assessment: the first SIC scalar, $\kappa_o^1 = \sum_{j \in B_o^{PTE}} \pi_j^*$.

### 4.2 Phase 2: Using the first scalar, $\kappa_o^1$

The solutions of TSP models (1) and (9) are the same because the STE1 program uses the SIC scalar, $\kappa_o^1$. The answers to the TGP models (2) and (10) differ because of the additional element, $\kappa_o^1 w_o^\#$. Comparing Eq. (3) with Eq. (12) would yield $(v_o^{PTE\#} x_o - u_o^{PTE\#} y_o) - (v_o^{STE1\#} x_o - u_o^{STE1\#} y_o) = \kappa_o^1 w_o^\#$.

Figure 1 in Section 5.1 depicts the 2D geometric intuition of the example's Step I virtual technologies of DMUs, called *interim virtual best technologies*. The locations of Kp and Tp, $(\alpha_o^\#, \beta_o^\#)$ and $(\hat{\alpha}_o^\#, \hat{\beta}_o^\#)$, are *DMU-o*'s current and goal *interim virtual values* in the PTE program, [see Eqs. (3) and (27)]. While in the STE1 program, locations of K1 and T1 are $(\alpha_{\omega o}^{1\#}, \beta_{\omega o}^{1\#})$ and $(\hat{\alpha}_{\omega o}^{1\#}, \hat{\beta}_{\omega o}^{1\#})$, [see Eq. (12) and (28)]. In the TSP models, the rectilinear distances between Kp and Tp equal $\delta_o^{PTE\#}$, and between K1 and T1 equal $\delta_o^{STE1\#}$. We have $\delta_o^{PTE\#} = \delta_o^{STE1\#} = (\sum_{i \in I} Q_{io}^\# + \sum_{r \in R} P_{ro}^\#) \tau_o^\#$, where $\tau_o^\# = \$1$. The line that links Kp and K1 is parallel to the *interim virtual price boundary*, the diagonal at the origin O. The location of the interim point of *DMU-o* in the TGP model of the STE1 program, IK1, is $(v_o^{STE1\#} x_o, u_o^{STE1\#} y_o)$. The rectilinear distance between IK1 and K1 equals the vScalar price, $\kappa_o^I w_o^\#$ $(= [1 - \gamma_o]\kappa_o^I w_o^\# + \gamma_o \kappa_o^I w_o^\#)$.

Step I of the PTE and STE1 programs have the relationships shown in Eq. (38). Note the best peers' intensities, $\pi_j^{PTE\#}$ and $\pi_j^{STE1\#}$, in the PTE and STE1 programs may be different; therefore, the targets for each performance index of the two programs are distant. The STE1 program comprehends the vScalar price to yield the other intensities.

$$\sum_{j \in B_o^{PTE}} \pi_j^\# = \sum_{j \in B_o^{STE1}} \pi_j^\# = \kappa_o^I, \Delta_o^{STE1\#} = \Delta_o^{PTE\#} = \delta_o^{STE1\#} = \delta_o^{PTE\#},$$

$$(Q_o^{PTE\#}, P_o^{PTE\#}) = (Q_o^{STE1\#}, P_o^{STE1\#}), \text{ and } (v_o^{PTE\#}, u_o^{PTE\#}) \neq (v_o^{STE1\#}, u_o^{STE1\#}). \tag{38}$$

Figure 2 depicts the 2D geometric intuition of virtual technology sets in the PTE and STE1 programs. The coordinate of *DMU-j* in the PTE program is $(\alpha_j^*, \beta_j^*)$, and in the STE1 program is $(\alpha_{\omega j}^{1*}, \beta_{\omega j}^{1*})$. Points Kp and Tp locate at $(\alpha_o^*, \beta_o^*)$ to $(\hat{\alpha}_o^*, \hat{\beta}_o^*)$, and points K1 and T1 locate at $(\alpha_{\omega o}^{1*}, \beta_{\omega o}^{1*})$ to $(\hat{\alpha}_{\omega o}^{1*}, \hat{\beta}_{\omega o}^{1*})$. In this particular example, we have $E_o^{STE1*} < E_o^{PTE*}$ because $w_o^* > 0$. When $w_o^* < 0$, we would have $E_o^{STE1*} > E_o^{PTE*}$. The goal relative efficiencies of the points T1 and Tp, $bE_o^{STE1*}$ and $bE_o^{PTE*}$, are equal to 1. The rectilinear distance between Kp and Tp equals $\Delta_o^{PTE*}$, and between K1 and T1 equals $\Delta_o^{STE1*}$. $(\sum_{i \in I} Q_{io}^* + \sum_{r \in R} P_{ro}^*)$ multiplies $\tau_o^{PTE*}$ equals $\Delta_o^{PTE*}$, and multiplies $\tau_o^{STE1*}$ equals $\Delta_o^{STE1*}$. The PTE and STE1 programs have the relationships shown in Eq. (39).

$$\tau_o^{STE1*} \neq \tau_o^{PTE*}, \Delta_o^{STE1*} \neq \Delta_o^{PTE*}, (v_o^{PTE*}, u_o^{PTE*}) \neq (v_o^{STE1*}, u_o^{STE1*}),$$

$$(Q_o^{PTE*}, P_o^{PTE*}) = (Q_o^{STE1*}, P_o^{STE1*}), \text{ and } 0 < E_o^{PTE*} \neq E_o^{STE1*} \leq 1. \tag{39}$$

Step II solutions of the VGP model in PTE and STE1 programs depend on $\tau_o^{PTE*}$ and $\tau_o^{STE1*}$, but the slack ratios are the same. $E_o^{PTE*}$ and $E_o^{STE1*}$ measured the relative efficiencies of *DMU-o* with the estimated slack ratios, reducing inputs and expanding outputs. $E_o^{PTE*}$ is the pure technical efficiency that ignores the effects in improving the scale of inputs and outputs, while $E_o^{STE1*}$ reflects the scale efficiency as shown in Eq. (34), $E_o^{STE1*} = T_o^{STE1*} - S_o^{STE1*}$. The characteristics of the PTE and STE1 programs clarified that scale efficiency is essential in evaluating *DMU-o*. Section 3.4 illustrated the derivations of scale efficiency [see Eqs. (34) and (35)].

One can observe the relationships between the PTE and STE1 programs in Figures 1 and 2. The PTE program provides the first scalar, $\kappa_o^1$. In Phase 3, the STE1 program will provide the second scalar, $\kappa_o^2$. In Phase 4, *DMU-o* will select the target scalar, $\hat{\kappa}_o$, between $\kappa_o^1$ and $\kappa_o^2$ to have achievable slacks of performance indices for implementation. The analysis of the PTE, STE1, and STE2 programs indirectly supports *DMU-o* for performance management in reality. The selected relative efficiency, $E_o^{STE1\wedge}$, has relationships with $E_o^{STE1*}$, $E_o^{STE2*}$, and $E_o^{PTE*}$, as depicted in Figures 2 and 3. One should not impose *DMU-o* to use the target scalar, $\hat{\kappa}_o$, equals 1, which is the convexity condition in VRS DEA programs.

## *4.3 Phase 3: Identifying the second scalar, $\kappa_o^2$*

Performing the standard sensitivity analysis of LP STE1 on the right-hand side, perturbation $\kappa_o^1$, would yield the allowable decrease and increase in $\kappa_o^1$, as shown in Eq. (40). Here, $\kappa_o^2$ denotes the second bound of the SIC scalar $\kappa_o^a$. The exact best peers in sets $B_o^{STE1}$ and $B_o^{STE2}$ are the basic variables of the two programs and will have different intensities, $\pi_j^* > 0$.

$$\kappa_o^2 = \kappa_o^1 - \text{(allowable decreasing of } \kappa_o^1\text{) or } \kappa_o^2 = \kappa_o^1 + \text{(allowable increasing of } \kappa_o^1\text{)}. \tag{40}$$

The STE1 and STE2 programs in Step I have solutions with the relationships summarized in Eq. (41).

$$B_o^{STE1} = B_o^{STE2}, (Q_o^{STE1*}, P_o^{STE1*}, \pi^{STE1*}) \neq (Q_o^{STE2*}, P_o^{STE2*}, \pi^{STE2*}), \tau_o^{STE1*} \neq \tau_o^{STE2*},$$

$$(v_o^{STE1\#}, u_o^{STE1\#}, w_o^{\#}) = (v_o^{STE2\#}, u_o^{STE2\#}, w_o^{\#}), \text{ and } (v_o^{STE1*}, u_o^{STE1*}, w_o^{*}) \neq (v_o^{STE2*}, u_o^{STE2*}, w_o^{*}). \tag{41}$$

Additionally, $\kappa_o^a$ represents an arbitrary scalar within the interval of the SIC scalars $\kappa_o^1 < \kappa_o^a < \kappa_o^2$ or $\kappa_o^1 > \kappa_o^a > \kappa_o^2$ when $w_o^* > 0$ and $w_o^* < 0$, respectively. $\Delta_{\omega o}^{STE2*} < \Delta_o^{PTE*}$ and $\Delta_{\omega o}^{STE2*} > \Delta_o^{PTE*}$ when $w_o^* > 0$ and $w_o^* < 0$, respectively. $\omega_o^{a*}$ is at the upper and lower bounds when $w_o^* < 0$ and $w_o^* > 0$, respectively. Phase 3 successfully identifies the exact second scale efficiency, $S_o^{STE2*}$. Figure 3 illustrates the solutions of the STE1 and STE2 programs of the numerical example in 2D geometric intuition. The analysts deal with the observed dataset in the first three phases. Phase 4 shall involve *DMU-o* in choosing the target SIC scalar between intervals of $\kappa_o^1$ and $\kappa_o^2$ to have an achievable improvement plan.

## 4.4 Phase 4: Select the target scalar $\hat{\kappa_o}$ within the interval of $[\kappa_o^1, \kappa_o^2]$

*DMU-o* should not blindly follow the mathematically analyzed best peers in the set of $B_o^{STE1}$, in which the number of peers will not be too many to confirm whether each peer and *DMU-o* have compatible production technologies. The ways to compare existing production technologies between *DMU-o* and peers are noteworthy and subject to *DMU-o*'s preference. *DMU-o* deletes incompatible peers from the observed dataset and recursively performs the first three phases.

Upon *DMU-o*'s confirmation of the best peers compatible with it, it could select the target scalar, $\hat{\kappa_o}$, within the interval of $\kappa_o^1$ and $\kappa_o^2$ and solve the STE^ program. That fulfills Criterion (v). The values $E_o^{STE2*}$, $T_o^{STE1*}$, $S_o^{STE1*}$, and $\Xi_o^{STE1*}$ could be at their upper and lower bounds if $w_o^* > 0$ and $w_o^* < 0$, respectively. Depending on the priorities of the assessment indices, *DMU-o* requires several trials to select a preferred target scalar, $\hat{\kappa_o}$, for the STE^ program. The target for each *input-i* and *output-r*, $\hat{x}_{io}$ and $\hat{y}_{ro}$, respectively, should be achievable for application. *DMU-o* will have the comprehensive assessment indices $E_o^{STE\wedge*}$, $T_o^{STE\wedge*}$, $S_o^{STE\wedge*}$, and $\Xi_o^{STE\wedge*}$. In Phase 4, *DMU-o* will obtain the preferred target scale efficiency, $S_o^{STE\wedge*}$. Sometimes, *DMU-o* is an efficient one, and it is the best performer.

Table 2 summarizes the numerical example solutions of the PTE, STE1, and STE2 programs, in which DMU-B and DMU-D are the best peers. The intensities, $\pi_B^*$ and $\pi_D^*$, are summarized in R3. R7 lists $\Xi_o^{STEa}$, $E_o^{STEa}$, $T_o^{STEa}$, and $S_o^{STEa}$, while R4 depicts the benchmark of each input and output. *DMU-o* exhibits decreasing and increasing bRTP when $w_o^* > 0$ and $w_o^* < 0$. In the best practice, *DMU-o* selects the upper bound of bRTP if each input and output benchmark is achievable. However, *DMU-o* could use Eq. (36) and (37) to analyze the interlinkage relationships between input and output indices for managing the performance indices.

## 5. Numerical examples and the 2D geometric intuitions

Table 1 shows the observed dataset of six DMUs (K, A, B, D, G, and H) with two inputs and outputs and with various measurement units. *DMU-o* is DMU-K and is evaluated according to the four-phase procedure involving the PTE, STE1, STE2, and STE^ programs. Table 2 summarizes the solutions of the VGA and CRS SBM programs. Because $w_o^\# > 0$, the SIC scalars' upper and lower bounds are $\kappa_o^1$ and $\kappa_o^2$. R8 and R9 list the virtual technologies of the DMUs in assessing *DMU-o*.

Table 1. The example data.

|  | $x_{1j}$ | $x_{2j}$ | $y_{1j}$ | $y_{2j}$ |
|---|---|---|---|---|
| *DMU-j* | ton | Hrs. | M$^3$ | °c |
| K | 1.6 | 145 | 1036 | 49 |
| A | 2.3 | 120  | 1327 | 97   |
| B | 1   | 29   | 567  | 89   |
| D | 1.9 | 281  | 2446 | 97   |
| G | 1.8 | 250  | 1794 | 57   |
| H | 2.5 | 100  | 1000 | 70   |

Table 2. *DMU-o* (=*DMU-K*) evaluations in the VGA and CRS SBM programs.

| Row | Solution's symbol | Unit | PTE-I | PTE-II | STE1-I | STE1-II | STE2-I | STE2-II | SBM | STE3-II |
|---|---|---|---|---|---|---|---|---|---|---|
| R1 | $\kappa_o^a$ | - | 1.5153 | 1.5153 | 1.5153 | 1.5153 | 0.5150 | 0.5150 |  | 1 |
|  | $\tau_o$ | $ | 1 | 0.179 | 1 | 0.256 | 1 | 0.5 | 0.531 | 0.338 |
| R2 | $\Delta_o, \delta_o$ | $ | 2.3010 | 0.4113 | 2.3010 | 0.5893 | 0.6643 | 0.3321 | 0.3893 | 0.4922 |
|  | $v_{1o}$ | $/ton | 2.8713 | 0.5133 | 0.6250 | 0.1601 | 0.6250 | 0.3125 | 0.6948 | 0.2110 |
|  | $v_{2o}$ | $/hrs | 0.0069 | 0.0012 | 0.0069 | 0.0018 | 0.0069 | 0.0034 | 0.0034 | 0.0023 |
|  | $u_{1o}$ | $/M^3 | 0.0022 | 0.0004 | 0.0011 | 0.0003 | 0.0011 | 0.0006 | 0.0008 | 0.0004 |
|  | $u_{2o}$ | $/°c | 0.0204 | 0.0036 | 0.0204 | 0.0052 | 0.0204 | 0.0102 | 0.0040 | 0.0069 |
|  | $w_o$ | $ | 0 | 0 | 1.6362 | 0.4190 | 1.6362 | 0.8181 | - | 0.5524 |
| R3 | $Q_{1o}$ | - | 0 | 0 | 0 | 0 | 0.4554 | 0.4554 | 0 | 0.2346 |
|  | $Q_{2o}$ | - | 0.5334 | 0.5334 | 0.5334 | 0.5334 | 0.2089 | 0.2089 | 0.5334 | 0.3662 |
|  | $P_{1o}$ | - | 0 | 0 | 0 | 0 | 0 | 0 | 0 | 0 |
|  | $P_{2o}$ | - | 1.7677 | 1.7677 | 1.7677 | 1.7677 | 0 | 0 | 1.7677 | 0.8571 |
|  | $\pi_B$ | - | 1.421 | 1.421 | 1.421 | 1.421 | 0.119 | 0.119 | 1.421 | 0.750 |
|  | $\pi_D$ | - | 0.094 | 0.094 | 0.094 | 0.094 | 0.396 | 0.396 | 0.094 | 0.250 |
| R4 | $\hat{x}_{1o}$ | ton | 1.6 | 1.6 | 1.6 | 1.6 | 0.8713 | 0.8713 | 1.6 | 1.225 |
|  | $\hat{x}_{2o}$ | hrs | 67.66 | 67.66 | 67.66 | 67.66 | 114.72 | 114.72 | 67.66 | 91.90 |
|  | $\hat{y}_{1o}$ | $M^3$ | 1036 | 1036 | 1036 | 1036 | 1036 | 1036 | 1036 | 1036 |
|  | $\hat{y}_{2o}$ | °c | 135.6 | 135.6 | 135.6 | 135.6 | 49.0 | 49.0 | 135.6 | 91.00 |
|  | $\gamma_o$ | - | 0.232 | 0.232 | 0.232 | 0.232 | 1.000 | 1.000 | 0.232 | 0.412 |
|  | $\omega_o^a$ | $ | 0 | 0 | 2.479 | 0.635 | 0.843 | 0.421 | 0 | 0.552 |
| R5 | $x_{1o}v_{1o}$ | $ | 4.594 | 0.821 | 1.000 | 0.256 | 1.578 | 0.789 | 1.112 | 0.553 |
|  | $x_{2o}v_{2o}$ | $ | 1.000 | 0.179 | 3.479 | 0.891 | 1.265 | 0.632 | 0.500 | 0.674 |
|  | $y_{1o}u_{1o}$ | $ | 2.293 | 0.410 | 1.178 | 0.302 | 1.178 | 0.589 | 0.806 | 0.398 |
|  | $y_{2o}u_{2o}$ | $ | 1.000 | 0.179 | 3.479 | 0.891 | 1.843 | 0.921 | 0.195 | 0.890 |
| R6 | $\alpha_o, \alpha_{\omega o}^a$ | $ | 5.594 | 1.000 | 3.905 | 1.000 | 2.000 | 1.000 | 1.612 | 1 |
|  | $\beta_o, \beta_{\omega o}^a$ | $ | 3.293 | 0.589 | 1.604 | 0.411 | 1.336 | 0.668 | 1.001 | 0.508 |
|  | $\hat{\alpha}_o, \hat{\alpha}_{\omega o}^a$ | $ | 5.061 | 0.905 | 3.371 | 0.863 | 1.336 | 0.668 | *1.470* | 0.797 |
|  | $\hat{\beta}_o, \hat{\beta}_{\omega o}^a$ | $ | 5.061 | 0.905 | 3.371 | 0.863 | 1.336 | 0.668 | *1.184* | 0.797 |
| R7 | $\Xi_o$ | - | 1.699 | 1.699 | 0.624 | 2.435 | 1.497 | 1.497 | 1.296 | 1.969 |
|  | $E_o$ | - | 0.589 | 0.589 | 1.604 | 0.411 | 0.668 | 0.668 | 0.621 | 0.508 |
|  | $T_o$ | - | 0.589 | 0.589 | 4.083 | 1.046 | 1.510 | 1.089 | 0.621 | 1.060 |
|  | $S_o$ | - | 0 | 0 | 2.479 | 0.635 | 0.843 | 0.421 | 0 | 0.552 |
| R8 | $\alpha_A, \alpha^a_{\omega A}$ | $ | 7.432 | 1.000 | 0.902 | 0.902 | 1.133 | 1.133 | 2.012 | 1.133 |
| | $\alpha_B, \alpha^a_{\omega B}$ | $ | 3.071 | 1.328 | 0.533 | 0.533 | 0.413 | 0.413 | 0.795 | 0.413 |
| | $\alpha_D, \alpha^a_{\omega D}$ | $ | 7.393 | 0.549 | 1.122 | 1.122 | 1.563 | 1.563 | 2.289 | 1.563 |
| | $\alpha_G, \alpha^a_{\omega G}$ | $ | 6.892 | 1.322 | 1.052 | 1.052 | 1.425 | 1.425 | 2.113 | 1.425 |
| | $\alpha_H, \alpha^a_{\omega H}$ | $ | 7.868 | 1.232 | 0.899 | 0.899 | 1.126 | 1.126 | 2.082 | 1.126 |
| R9 | $\beta_A, \beta^a_{\omega A}$ | $ | 4.917 | 0.879 | 0.796 | 0.796 | 1.853 | 0.926 | 1.418 | 0.926 |
| | $\beta_B, \beta^a_{\omega B}$ | $ | 3.071 | 0.549 | 0.533 | 0.533 | 0.825 | 0.413 | 0.795 | 0.413 |
| | $\beta_D, \beta^a_{\omega D}$ | $ | 7.393 | 1.322 | 1.122 | 1.122 | 3.125 | 1.563 | 2.289 | 1.563 |
| | $\beta_G, \beta^a_{\omega G}$ | $ | 5.134 | 0.918 | 0.723 | 0.723 | 1.568 | 0.784 | 1.623 | 0.784 |
| | $\beta_H, \beta^a_{\omega H}$ | $ | 3.642 | 0.651 | 0.560 | 0.560 | 0.930 | 0.465 | 1.056 | 0.465 |

Sections 3.4 and 4.2 describe the relationships of $E_o^{PTE*}$, $E_o^{STE1*}$, $T_o^{STE1*}$, and $S_o^{STE1*}$. Figure 1 depicts the geometric intuition of the PTE and STE1 programs. The best peers B and D are located on the interim virtual price boundary, the diagonal at the origin O with a slope of 1. The rectilinear distance between K1 and T1 is the total VGP in the STE1 model, $\delta_o^{STE1\#}$, which is the sum of $(Q_{1o}^* + Q_{2o}^*)\tau_o^{STE1\#}$ and $(P_{1o}^* + P_{2o}^*)\tau_o^{STE1\#}$.

In Figure 2, the virtual technologies of each *DMU-j* in the STE1 and PTE programs are points, with *1* and *p* appended to the labels, $(\alpha_j^*, \beta_j^*)$ and $(\alpha_{\omega j}^{1*}, \beta_{\omega j}^{1*})$, respectively. Points Kp and K1 are on the vertical normalization line (N, E). The targets Tp and T1 are located on the virtual best-efficiency boundary at $(\hat{\alpha}_o^*, \hat{\beta}_o^*)$ and $(\hat{\alpha}_{\omega o}^{1*}, \hat{\beta}_{\omega o}^{1*})$, respectively. The total virtual slack price of the inputs and outputs in the STE1 model is the rectilinear distance between points K1 and T1, $\sum_{i \in I} Q_{io}^* \tau_o^*$ plus $\sum_{r \in R} P_{ro}^* \tau_o^*$.

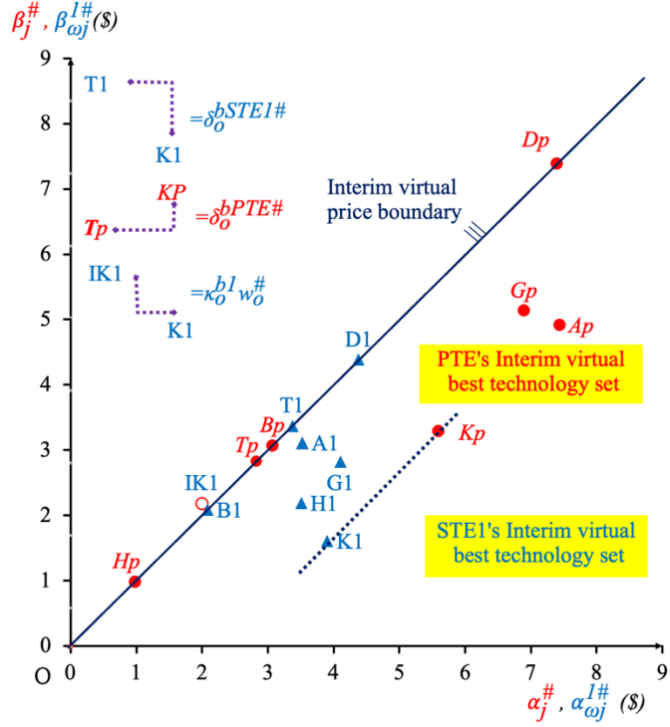

Figure 1: The PTE and STE1's Step 1 interim solutions.

In Eq. (25), when $DMU\text{-}j$ has $\Delta_j^{STEa*} > 0$ because ($\alpha_{\omega j}^{a*} > 0$, $\beta_{\omega j}^{a*} < 0$, and $\beta_{\omega j}^{a*}/\alpha_{\omega j}^{a*} < 0$) and ($\alpha_{\omega j}^{a*} < 0$, $\beta_{\omega j}^{a*} < 0$, and $\beta_{\omega j}^{a*}/\alpha_{\omega j}^{a*} > 0$), $DMU\text{-}j$ is located in the fourth and third quadrants, respectively. Therefore, eliminating the inefficient $DMU\text{-}j$ from $(X, Y)$ does not affect the evaluation. In this case, $DMU\text{-}o$ exhibits decreasing bRTP because $w_o^* > 0$. $DMU\text{-}o$ selects the STE1 model and $\hat{\kappa}_o = \kappa_o^I$ if the benchmark of each input and output is achievable in practice.

In Figure 3, for simplicity, the three vectors $\overrightarrow{O,K1}$, $\overrightarrow{AP1,O}$, and $\overrightarrow{AP1,K1}$ are found in Figure 2. The anchor point AP2 is located at $[(1-\gamma_o)\omega_o^{2*}, -\gamma_o\omega_o^{2*}]$. Like Eq. (35), the sum of the three graphically geometrical vectors, $\overrightarrow{O,K2} - \overrightarrow{AP2,O} = \overrightarrow{AP2,K2}$, is equivalent to $\Delta_o^{STEa*} - \omega_o^{a*} = (v_o^* x_o - u_o^* y_o)$ [see Eq. (31)]. Because $w_o^\# > 0$ or $\kappa_o^a w_o^\# > 0$, Table 2 presents the decreasing order of $S_o^{STE1*} < S_o^{STE2*}$ and the increasing order of $E_o^{STE1*} > E_o^{STE2*}$ and $\Xi_o^{STE1*} > \Xi_o^{STE2*}$ with an increasing $\kappa_o^a$ (see R7 in Table 2). Because $w_o^* > 0$, the anchor points AP1 and AP2 are located at $[(1-\gamma_o)\omega_o^{1*}, -\gamma_o\omega_o^{1*}]$ and $[(1-\gamma_o)\omega_o^{2*}, -\gamma_o\omega_o^{2*}]$ in the fourth quadrant, respectively.

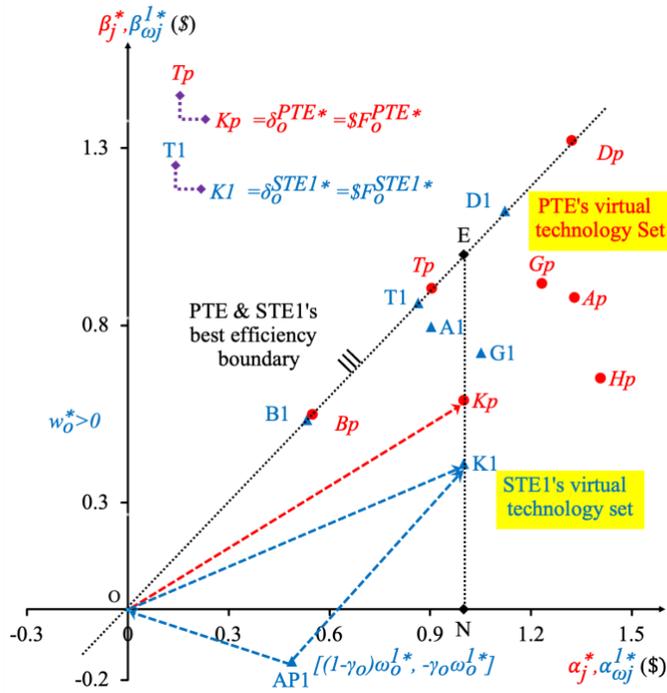

Figure 2: Final solutions of the PTE and STE1 programs.

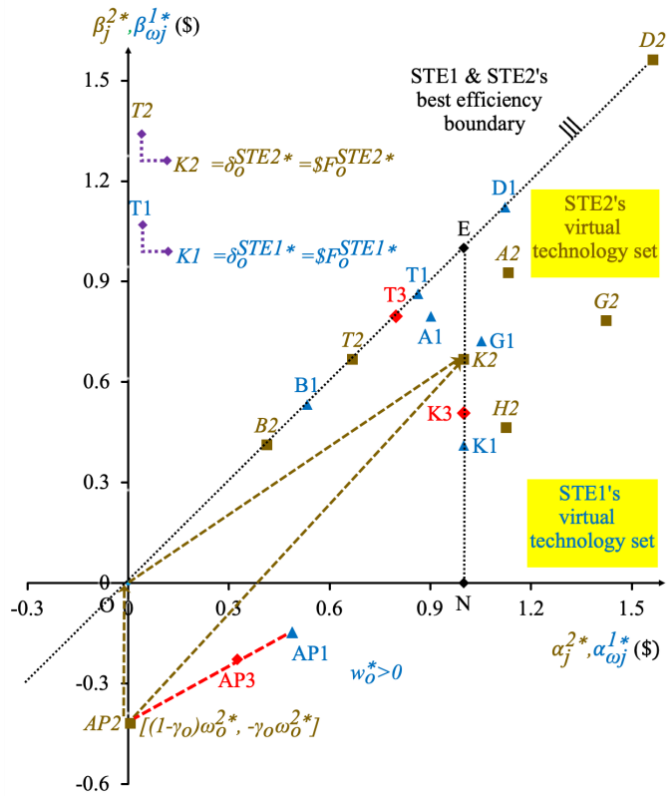

Figure 3: The STE1 and STE2's final solutions.

Using the target SIC scalar $\hat{\kappa}_o$ in the STE^ model, points AP^, K^, and T^ will be located on the line segment between points AP1 and AP2, K1 and K2, and T1 and T2, respectively. The three vectors, $\overrightarrow{O,K\hat{}}$, $\overrightarrow{AP\hat{},O}$, and $\overrightarrow{AP\hat{},K\hat{}}$, shall be similar to Eq. (35). For instance, $\kappa_o^3 = 1$ in the STE3 model, the alternative $\kappa_o^3$ is feasible because $\kappa_o^1 > \kappa_o^3 > \kappa_o^2$. If $\kappa_o^3$ is not within the interval of $\kappa_o^1$ and $\kappa_o^2$, the alternative is not feasible. The right-hand column of Table 2 lists the solutions of the CRS SBM model, and the solutions in R3 and R4 are the same as in the PTE-II model. In R2, $E_o^{SBM} = 0.3893$, $v_o^* x_o = 1.612 > 1$, and the relative efficiency in R5, $E_o^{Rel} = u_o^* y_o / v_o^* x_o = 0.621$, are seen. Because $E_o^{SBM} < E_o^{Rel}$, the CRS SBM model has incomplete solutions. The exact relative efficiency is $E_o^{PTE} = 0.589$ (see R7). In R6, the goal relative efficiency, $\hat{\beta}_o/\hat{\alpha}_o$ (= 1.184/1.47 ≠ 1), is not equal to 1, which indicates that *DMU-o* could not hold Criterion 2. Thus, the model has an incomplete solution.

In contrast, the PTE and STE3 programs provide exact solutions for the CRS and VRS NONR DEA programs.

## 6. Discussion

Inspired by the CRS and VRS NONR DEA programs, we developed the PTE and STEa programs that satisfy Criterion (i)–Criterion (iv). In the four-phase VGA procedure, the STE^ program fulfills Criterion (v) for *DMU-o*. This research aims to identify the best benchmark for each input and output for implementation.

Based on the Simplex method (Dantzig et al., 1997), VGA programs have no assumptions and particular theorems. One could perform the post-optimality analysis in the STE^ program with the fundamental knowledge of linear programming. The changes in the coefficients of the program affect the optimal solution. Such analysis is an essential part of evaluating *DMU-o* in practice.

Ranking DMUs according to the STE^ solutions of each *DMU-o* may not be adequate because they subjectively select the $\hat{\kappa}_o$ and have different best peers. Currently, we are not ranking DMUs with relative efficiency equal to 1. In the worst practice assessment, DMUs have other rankings.

Generally, $Q_o^*$ is between 0 and 1, and $P_o^*$ may be greater than 1. Sometimes one may add the condition $P_{ro} \leq 1$, $r = 1, \ldots, s$. *DMU-o* may also have additional constraints with lower and upper bounds to the slack ratios, $\underline{Q}_o \leq Q_o \leq \overline{Q}_o$ and $\underline{P}_o \leq P_o \leq \overline{P}_o$, respectively, to meet its requirements (Liu & Shih, 2022a).

To effectively use the VGA procedure, fundamental performance indices must be identified. Limited observation bias occurs in VGA; however, this bias can be reduced if enough of the highest or near-highest performers are included. Constraints may be added to meet the practical requirements based on the VGA programs. Several practical problems have been discussed in the DEA literature, such as the Malmquist Index, network DEA, free-disposal hull, imprecise data, data mixed with 0s and negatives, input/output sharing, and input/output selection (Panwar et al., 2022). Employing the PTE and STEa programs could solve the problems of multi-criteria decision-making methods [see Stojčić et al. (2019)].

Shih (2022) described the four-phase VGA procedure used to evaluate performance in the worst practice. Liu and Shih (2022b) used the four-phase VGA procedure in the most and least favorable practices to estimate *super-efficiency* and *hypo-inefficiency*.

**Acknowledgments:** This study was inspired by the DEA, SFA, DDF, and multi-criteria decision-making methods (MCDM). We would like to thank the referees for their valuable comments.

**Funding:** This research received no specific grant from public, commercial, or not-for-profit funding agencies.

**Disclosure:** We have no competing interests to declare.

**Data availability statement:** We confirm that the data supporting the findings of this study are available within the article.